\documentclass[reqno]{amsart}
\usepackage{amsfonts,amsthm,amsmath,amssymb,amscd,amsmath,epsf,latexsym,enumerate,mathrsfs}
\usepackage{graphicx}
\usepackage[utf8x]{inputenc}
\newtheorem{theorem}{Theorem}

\newtheorem{lemma}[theorem]{Lemma}

\newtheorem{problem}{Problem}

\newcommand{\ee}{\end{equation}}

\begin{document}
          \numberwithin{equation}{section}

          \title[reality of zeros for polynomials obeying a three-term recurrence]
          {CRITERION OF REALITY OF ZEROS IN A POLYNOMIAL SEQUENCE SATISFYING A THREE-TERM RECURRENCE RELATION}

\author[I.~Ndikubwayo]{Innocent Ndikubwayo}
\address{Department of Mathematics, Stockholm University, SE-106 91
Stockholm,         Sweden}

\email {innocent@math.su.se}

\address{Department of Mathematics, Makerere University, 7062 Kampala, Uganda}

\email {ndikubwayo@cns.mak.ac.ug}


\begin{abstract} This paper establishes the necessary and sufficient conditions for the reality of all the zeros in a polynomial sequence  $\{P_i\}_{i=1}^{\infty}$ generated by a three-term recurrence relation
$P_i(x)+ Q_1(x)P_{i-1}(x) +Q_2(x) P_{i-2}(x)=0$ with the standard initial conditions $P_{0}(x)=1, P_{-1}(x)=0,$ where $Q_1(x)$ and $Q_2(x)$ are arbitrary real  polynomials.

\end{abstract}

\maketitle

\section{Introduction}

 Asymptotic root distributions for sequences of univariate polynomials has been a topic of study in analysis for many decades \cite{RS}. In particular, sequences of polynomials with all real zeros are important  in many branches of mathematics. Such polynomials possess several nice properties. For example, if a polynomial $P(x)=\sum_{i=0}^{n}b_ix^i$ is real-rooted and has nonnegative coefficients, then the sequence $\{b_i \}_{i=0}^n$ is log-concave, i.e $b_i^2 \geq  b_{i+1}b_{i-1}$ for all $1 \leq j< n$ \cite{PB}.  This log-concavity implies that $\{b_i \}_{i=0}^n$ is  unimodal, whereby the sequence increases to a greatest value (or possibly two consecutive equal values) and then decreases \cite{PB}.  In addition, polynomials with real zeros  are closed with respect to differentiation and the zeros of derivative interlace with the zeros of the polynomial. A good source of information about real-rooted polynomials is a book V.P. Kostov, ``Topics on hyperbolic polynomials in one variable" \cite{BOo}.

\medskip 
In this paper, we discuss some cases of the problem under which conditions polynomials satisfying a finite  linear recurrence relation  have all real roots. 
Our general set-up is as follows. Fix complex-valued polynomials $Q_1(x), \ldots, Q_k(x)$ and consider a finite linear recurrence relation of the form
\begin{equation} \label{eqn1}
P_i(x)+ Q_1(x)P_{i-1}(x)+ Q_2(x)P_{i-2}(x)+ \dots + Q_k(x)P_{i-k}(x)=0,~~~~~~~i=1,2,\dots
\end{equation}
 with the standard initial conditions 
 
 \begin{equation}\label{eqn2}
 P_{0}(x)=1, P_{-1}(x)= P_{-2}(x)=  \dots =  P_{-k+1}(x)=0.
 \end{equation}

The generating function for the polynomial sequence $\{P_i\}_{i=1}^{\infty}$ of this recurrence is given by

\begin{equation*}
 \sum_{i=0}^{\infty}P_i(x)t^i= \frac{1}{1+Q_1(x)t+ Q_2(x)t^2+ \dots + Q_k(x)t^k}.
 \end{equation*}\\

One of the well-known results in this area is a description of the accumulation set for the zeros of $P_i(x)$ when $i \to \infty$ provided by the theorem by Beraha, Kahane and Weiss, see \cite{KK}.
It asserts that the support of the limiting root-counting measure coincides with  the  following set.
Let $Q_1,\dots, Q_k$  be complex polynomials as given in Equation (\ref{eqn1}) above. Define a curve $\Gamma_Q \subset \mathbb{C}$ consisting of all values of $x$ such that the characteristic equation 
\begin{equation} \label{eqn3}
 1+Q_1(x)t+ Q_2(x)t^2+ \dots + Q_k(x)t^k=0
 \end{equation}
 has at least two roots $t_1, t_2$ for which

\medskip
\noindent
(a) $|t_1|=|t_2|$;

\medskip
\noindent
(b) $|t_1|$ is the minimum among the absolute values of all roots.
\medskip

\begin{theorem} \cite{KK} Suppose that $ \{P_i(x)\}$ satisfies (\ref{eqn1}), (\ref{eqn2}) and (\ref{eqn3}). 
Suppose further that $\{P_i(x)\}$ satisfies  no recursion of order less than $k$ and that there is no constant $\omega \in \mathbb{C}$ of unit modulus for which $ t_r = \omega t_s $ for some $r \neq s.$ Then the zeros of $P_i(x)$ accumulate along
 the curve $\Gamma_Q$ as $i \to \infty$.
\end{theorem} 

This result provides a description of the asymptotic behaviour of the roots of $P_i(x)$. However recently Tran \cite{T} has found a number of cases when the zeros of $P_i(x)$ actually lie on the limiting curve $\Gamma_Q$ for all or for all sufficiently large $i$. In particular, he has proven the following results.

\begin{theorem} \cite{T}
Let $ \{P_i(x)\}$ be a polynomial sequence whose generating function is 
\begin{equation*}
 \sum_{i=0}^{\infty}P_i(x)t^i= \frac{1}{1+Q_1(x)t+ Q_2(x)t^2},
 \end{equation*}
 where $Q_1(x)$ and $Q_2(x)$ are polynomials in $x$ with complex coefficients. All the zeros of every polynomial in the sequence $\{P_i(x)\}$ which satisfy  $Q_2(x) \neq 0$ lie on the curve $\Gamma_Q$  defined by  \end{theorem}
 \begin{equation} \label{im}
 \mbox{Im} \left( \frac{ Q_1^2(x)}{Q_2(x)}\right)=0 ~~~~\mbox{and}~~~~ 0 \leq \mbox{Re}\left( \frac{ Q_1^2(x)}{Q_2(x)}\right)\leq 4.  \end{equation}

$ Moreover,~these~~ zeros~~ become~~ dense~~ in~~  \Gamma_Q~~ when~~ i \to \infty.$

Theorem $2$ covers the special case of (\ref{eqn1}) and (\ref{eqn2}) for the polynomials generated by the recurrence 
\begin{equation*}
P_i(x)+ Q_1(x)P_{i-1}(x)+ Q_2(x)P_{i-2}(x)=0 \end{equation*}
with the standard initial conditions $P_0(x)=1$ and $P_{-1}(x)=0$. 

A more general result of Tran is as follows.
\begin{theorem} \label{thm:three}\cite{TI}
Let $\{P_i(x)\}$ be a polynomial sequence with generating function 
\begin{equation*}
 \sum_{i=0}^{\infty}P_i(x)t^i= \frac{1}{1+Q_1(x)t+ Q_2(x)t^k},
 \end{equation*}
 where $Q_1(x)$ and $Q_2(x)$ are polynomials in $x$ with complex coefficients. Then there exists a constant $C=C(k)$ such that for all $i > C$, all the zeros of $P_i(x)$ which satisfy $Q_2(x) \neq 0$ lie on the curve $\Gamma_Q$  defined by  \end{theorem}
 \begin{equation*}
 \mbox{Im}\left( \frac{ Q_1^k(x)}{Q_2(x)}\right)=0 ~~~~\mbox{and}~~~~ 0 \leq (-1)^n\mbox{Re}\left( \frac{ Q_1^k(x)}{Q_2(x)}\right)\leq \frac{k^k}{(k-1)^{k-1}}.  \end{equation*}
$ Moreover,~these~~ zeros~~ become~~ dense~~ in~~  \Gamma_Q~~ when~~ i \to \infty.$

\medskip
In the present paper,  we want to characterise a situation as above that gives rise to polynomial sequences with only real zeros. 

\begin{problem}\label{prob:main}
{\rm} In the above notation, consider the recurrence relation \begin{equation*}P_i(x)+ Q_1(x)P_{i-1}(x)+ Q_2(x)P_{i-2}(x)=0 \end{equation*} with the standard initial conditions, \begin{equation*} P_{0}(x)=1, P_{-1}(x)=0\end{equation*} where $Q_1(x)$ and $Q_2(x)$ are arbitrary real polynomials. Give necessary and sufficient conditions on $(Q_1(x), Q_2(x))$  guaranteeing that all the zeros of $P_{i}(x)$ will be real for all $i$. 
\end{problem}

\medskip\noindent
To formulate our main result, we need to look at the curve defined by the first condition in $(\ref{im}).$ We shall view  $\mathbb{C}P^1$ as $\mathbb{C}\cup \{\infty\}$, the extended complex plane and $\mathbb{R}P^1$ as  the union of the real line in $\mathbb{C}$ with $\{\infty\}$.
\medskip
Let $f:\mathbb{C}P^1 \to \mathbb{C}P^1$ be the rational function defined by $f=\frac{Q_1^2(x)}{Q_2(x)}$ 
where $Q_1(x)$ and $Q_2(x)$ are real polynomials. Denote by $\widetilde{\Gamma_Q} \subset \mathbb{C}P^1$ the curve given by Im$(f)=0$, that is $$\widetilde{\Gamma_Q}=f^{-1}(\mathbb{R}P^1).$$ 
We note that for real polynomials $Q_1(x)$ and $Q_2(x)$, the curve $\widetilde{\Gamma_Q}$ contains $\Gamma_Q$ since $[0,4] \subset \mathbb{R}P^1.$

\medskip

\begin{lemma} \label{four}
The curve $\widetilde{\Gamma_Q}$ has the following properties:\\

\noindent
{\rm (a)} $\widetilde{\Gamma_Q} \supset \mathbb{R}P^1;   $\\

\noindent
{\rm (b)} $\widetilde{\Gamma_Q}$ is invariant under complex conjugation;\\

\noindent
{\rm (c)} except $\mathbb{R}P^1$, $\widetilde{\Gamma_Q}$ might contain ovals disjoint from $\mathbb{R}P^1$ (which come in complex-conjugate pairs) and ovals crossing $\mathbb{R}P^1$ which are mapped to themselves by complex conjugation;\\

\noindent
{\rm (d)} the intersection points of the second type of ovals with $\mathbb{R}P^1$ are exactly the real critical points of $f$. \\
\end{lemma}
Figures $1$ and $2$  illustrate the properties of $\widetilde{\Gamma_Q}$ claimed in  Lemma 4 $(a)-(d)$.  The main result of the present paper is as follows.

\begin{figure}[h]\begin{center}
$ 
\begin{array}{c}
\includegraphics[height=6cm, width=6cm]{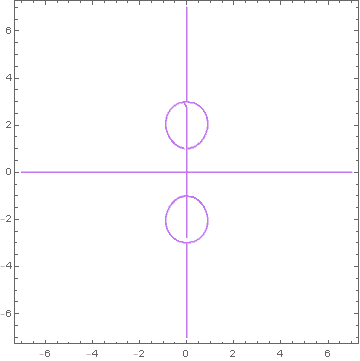}
\end{array}
$
\end{center}
\caption{The curve $\widetilde{\Gamma_Q}$  for  $f=\frac{Q_1^2(x)}{Q_2(x)}$ where $Q_1(x) = x^2 + 1$ and 
$Q_2(x) = x^2 +6$.}
\end{figure}

\begin{figure}[]\begin{center}
$ 
\begin{array}{c}
\includegraphics[height=6cm, width=6cm]{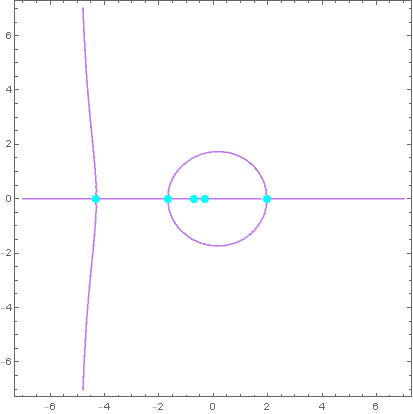}
\end{array}
$
\end{center}
\caption{The curve $\widetilde{\Gamma_Q}$  for  $f=\frac{Q_1^2(x)}{Q_2(x)}$ where   $Q_1(x) = x^2 +5x+3$ and $Q_2(x) = 5x^2-1.$ The green points are the real critical points of $f$.}

\end{figure}

\begin{theorem}\label{th:main} Let $\{P_i(x)\}$ be a sequence of polynomials whose generating function is 

\begin{equation}
 \sum_{i=0}^{\infty}P_i(x)t^i= \frac{1}{1+Q_1(x)t+ Q_2(x)t^2},
 \end{equation}
 where $Q_1(x)$ and $Q_2(x)$ are arbitrary coprime real polynomials in $x$. Then, for all positive integers $i$, all the zeros of $P_{i}(x)$ are real if and only if  the following conditions are satisfied:\\

{\rm (a)} The polynomial $Q_1(x)$ must have all real and simple zeros. \\

{\rm (b)} No ovals $\gamma$ of $\widetilde{\Gamma_Q}$ disjoint from $\mathbb{R}P^1$ should exist.\\

{\rm (c)} All the zeros of the discriminant $D(x)$ of the characteristic polynomial $1 + Q_1(x)t + Q_2(x)t^2$ must be real.\\

{\rm (d)} No real critical values of $f$ should belong to the interval $(0,4)$.\\

{\rm (e)} The polynomial $Q_2(x)$ must be non-negative at the zeros of $Q_1(x).$ 
\end{theorem}
\medskip
\noindent
{\bf Remark.} The situation when $Q_1(x)$ and $Q_2(x)$ have a common real zero is not interesting to consider since from  $P_i(x)+ Q_1(x)P_{i-1}(x)+ Q_2(x)P_{i-2}(x)=0$, such a zero would necessarily be a zero of $P_i(x)$ for all $i$. \\

{\bf Acknowledgements.}  I am sincerely grateful to my advisor Professor Boris Shapiro who introduced me to the problem and for the many discussions and guidance of the above problem. I also thank  Professor Rikard B\o gvad and Dr. Alex Samuel Bamunoba for the discussions and guidance. I acknowledge and appreciate the financial and academic support from Sida Phase-IV bilateral program with Makerere University 2015-2020 under project 316 `Capacity building in mathematics and its applications'.

\section{PROOFS}

Let us begin with the following definitions and remarks.
\medskip

\noindent
\textbf{Definition 1.} For a non-constant rational function $R(x) = \frac{P(x)}{Q(x)}$, where $P(x)$ and $Q(x)$ are coprime polynomials, the degree of $R(x)$ is defined as the maximum of the degrees of $P(x)$ and  $Q(x)$.
\\

\textit{Equivalently,} the degree of $R(x)$ is the  number of distinct preimages of any generic point. \\

\noindent
\textbf{Definition 2.} A  point $x_0$ is called a critical point of $R(x)$ if $R(x)$ fails to be injective in a neighbourhood of $x_0$, that is,  $R'(x_0)=0$. A critical value of $R(x)$ is the image of a critical point. The order of a critical point $x_0$ of $R(x)$ is the order of the zero of $R'(x)$
at $x_0$.\\

Now, if $d$ is the degree of $R(x)$ and if $w$ is not a critical value then $R^{-1}(w)=\{x_1, x_2, \dots, x_d\}$ with $x_i \neq x_j$ for all $i \neq j$. Since the points $x_j$ are non-critical there is a neighbourhood of each of these points such that $R(x)$ is injective on that neigbourhood. The function $R(x)$ has $2d-2$ critical points  in $\mathbb{C}P^1$ counting multiplicities. This follows from the fact that in the complex plane, deg $(R')=$ deg$( P'Q-Q'P )=$ deg$(P) +$ deg$(Q)-1$ while the order of the critical point at infinity is $|$deg$(P) -$ deg$(Q)|-1,$  \cite{GM}. \\

\noindent
   \textbf{Definition 3.} Given a pair $(P, Q)$ of polynomials, define their Wronskian as $$W(P, Q) = P'Q-Q'P.$$  One interesting thing about the Wronskian is that if $P$ and $Q$ are coprime, then the zeros of $W(P, Q)$ are exactly the critical points of the rational map $R=P/Q$.
In \cite{BO} we find that if $P$ and $Q$ have all real, simple and interlacing zeros,  then all zeros of
$W(P, Q)$ are non-real. In addition, if we know that $\alpha$ is a zero of $R$ of multiplicity $\geq 2$, then $\alpha$ is also a (multiple) zero of the Wronskian. More information about the Wronski map can be found in \cite{KD,BO}. \\

\medskip

\begin{proof}[Proof of Theorem 5]
(a) Substitution of the initial conditions $P_0(x) =1$ and $P_{-1}(x) =0$ in the recurrence relation
$$P_i(x)+ Q_1(x)P_{i-1}(x) +Q_2(x) P_{i-2}(x)=0$$ gives for $i=1$, 
$$P_1(x)+ Q_1(x)P_{0}(x) +Q_2(x) P_{-1}(x)=0$$  or $$P_1(x)=-Q_1(x).$$ Therefore $Q_1(x)$ must have all its zeros real since we require all the zeros of $P_i(x)$ to be real for all $i$ and in particular for $i=1$. 
These zeros must be simple (see part $(e)$ for the justification).
\medskip

\noindent
(b) Suppose there exists an oval $\gamma$ of $\widetilde{\Gamma_Q}$ which does not intersect $\mathbb{R}P^1$. From Lemma \ref{four} $(c)$, $\gamma$ is the type one oval contained in $\widetilde{\Gamma_Q}.$ We note that all the points on $\gamma$ and its interior are of the form $z=x+iy$ where $x, y \in \mathbb{R}, y \neq 0$ and this is a connected component  with $\gamma$ as its boundary. This component is mapped by $f$ onto the half plane with degree $\geq 1$ depending on the number of critical points of $f$ it strictly contains. The boundary $\gamma$ of the component, is mapped onto $\mathbb{R}P^1$ (the boundary of the half plane). In particular, the image $f(\gamma)$ covers the interval
$[0,4] \subset \mathbb{R}P^1$. Therefore $\Gamma_Q=f^{-1}([0,4])$ must contain an arc of the boundary $\gamma$.  From Theorem 2 \cite{T}, all the zeros of $P_i(x)$ are contained in $\Gamma_Q$  for all $i$ and are dense in $\Gamma_Q$. Now since we require all the zeros of $P_{i}(x)$ to be real, it must be that $\Gamma_Q \subseteq \mathbb{R}P^1$. This is not possible as we already have that $\Gamma_Q$ contains an arc  of $\gamma$ yet this arc is not contained in $\mathbb{R}P^1$, hence a contradiction.
\medskip

\noindent
(c) It is known \cite{Bg} that the endpoints of $\Gamma_Q$ are the points where $ t_1=t_2$ (see the notation in Theorem 2). In our case equation (\ref{eqn3}) has degree $2$ in $t$. Therefore every $x$, for which the roots of (\ref{eqn3}) coincide belongs to $\Gamma_Q$. These $x$ are exactly the zeros of $D(x)=Q_1^2(x)-4Q_2(x).$ Since we require that  $\Gamma_Q \subset \mathbb{R}P^1 $, all the zeros of $D(x)$ must be real.
\medskip

\noindent
(d) Suppose there exists a critical value $w \in (0,4)$. Then there must exist a real critical point $x_c$ such that $f(x_c)=w.$ Clearly, $x_c$ is a point on $\Gamma_Q$. 
It is known \cite{GM} that a point $x_c \in \mathbb{C}P^1 $ is a critical point of order $k$ for a rational function $R(x)$ if and only if there are open sets $U$ containing $x_c$ and $V$ containing $w=R(x_c)$ such that each $w_0 \in V, w \neq w_0$ has exactly $k+1$ distinct preimages in $U$. 

 In our scenario let $x_c$ be such a critical point of order $k$ for $f(x)$ and $V$ be the real interval  $(w -\epsilon, w +\epsilon)$  for sufficiently small $\epsilon >0$. Then $V \subset (0,4)$. Note that since $x_c$ is a critical point of order $k$ for $f(x)$  and any point  $z \in V, z \neq w$ has exactly $k+1$ distinct preimages in $U$, we have at $x_c$, locally  $U= f^{-1}(V)$ the preimage of $V$ consists of $k+1$ distinct curves (arcs) with a common intersection only at $x_c$. One of these curves is a line segment on the real line while the remaining $k$ curves are  complex i.e apart from $x_c$, points on these $k$ curves are of the form  $z=x+iy$ where $x, y \in \mathbb{R}, y \neq 0$.
 
Now since complex arcs are formed in the preimage of $V$, then some of the zeros of $P_i(x)$ will be contained in the complex arcs since all the zeros of $P_i(x)$ are contained in ${\Gamma_Q}=f^{-1}([0,4] \supset V)$ and are dense there as $i \to \infty$. This contradicts our requirement that ${\Gamma_Q}$ is contained in $\mathbb{R}P^1 $. Therefore in order to have all the real roots of $P_i(x)$ for all $i$, no real critical values of $f$ can be in the real interval $(0,4)$. Otherwise the condition that ${\Gamma_Q} \subset \mathbb{R}P^1$ cannot hold.

\medskip

\noindent
(e) Let $x_0$ be a zero of $Q_1(x)$, i.e $Q_1(x_0)=0$. Note that $Q_2(x)$ and $Q_1(x)$ do not have a common zero since they are coprime. Therefore at the point $x_0$, we have $Q_2(x_0) \neq 0.$ It remains to show that $Q_2(x_0)>0$. 
 We note that all the zeros of $Q_1(x)$ are critical points of $f$. In addition, they belong to $\Gamma_Q$ because at $x_0$, we have that $f(x_0) =Q_1^2(x_0)/Q_2(x_0)= 0$;  therefore both the real and the imaginary part of $f$ vanish, hence satisfying (\ref{im}) of Theorem $2$. 
Suppose that $x_0$ is a simple critical point of $f$ and let  $Q_2 > 0$ at $x_0$. Then locally $ f= Q_1^2/Q_2 \geq 0$ on the interval in $\mathbb{R}$ and if $Q_2 < 0$ at $x_0$, then locally  $f= Q_1^2/Q_2 \geq 0$ on the complex arc, i.e there exists an interval $I \subset [0,4]$ such that  $f^{-1}(I)$ contains a complex arc. 
On the other hand, if $x_0$ is  a critical point of order $> 1$, then  locally $ f= Q_1^2/Q_2 \geq 0$ on some complex arc irrespective of whether $Q_2>0$ or $Q_2<0$. However, it is known that the zeros of $P_i(x)$ are contained in $\Gamma_Q=f^{-1}[0, 4]$ and are dense there as $i \to \infty$, and so for the reality of all the zeros of $P_i(x)$ we require that $\Gamma_Q \subset \mathbb{R}P^1$.  This is not possible if $Q_2<0$ at simple zeros of $Q_1(x)$ or when zeros of $Q_1$ have multiplicity $> 1$ since in either case there will be some zeros of $P_i(x)$ on the complex arc which is a contradiction. Therefore the polynomial $Q_2(x)$ must be non-negative at the  zeros of $Q_1(x)$ as a necessary condition for the reality of all the zeros of $P_i(x)$ for all $i$. Furthermore all the zeros of $Q_1(x)$ must be simple otherwise as explained above some zeros of $P_i(x)$ would be on the complex arc. (This last part settles part $(a)$ of the Theorem where we require all the zeros of $P_i(x)$ to be simple).
\end{proof}

\noindent
\textbf{Remark.} Each of the conditions of Threorem 5 $(a)-(e)$ is only a necessary (and not a sufficient) condition for the reality of all the zeros of $P_i(x)$. To guarantee reality of all the zeros of $P_i(x)$ for all $i$, all the five conditions must be satisfied simultaneously.  Below are some of the examples illustrating this claim. \\

\medskip\noindent
\textbf{Example 1.} Consider the sequence of polynomials $\{P_{i}(x)\}$ generated by the rational function

\begin{equation*}
 \sum_{i=0}^{\infty}P_i(x)t^i= \frac{1}{1+(-x^2+2x)t+ (5x^2-1)t^2}.
 \end{equation*}
  The corresponding $f$ is given by $$f(x)=\frac{(-x^2+2x)^2}{5x^2-1}.$$ The zeros of $Q_1(x)$ are $0$ and $2$ which are real and simple (Theorem $5 (a)$). There are no ovals disjoint from $\mathbb{R}P^1$ (Theorem $5 (b)$). The discriminant $D(x)= x^4-4x^3 - 16x^2+4$ has only real zeros (Theorem $5 (c)$). These are $-2.39337, -0.54374, 0.47570$ and $6.46141$ (corrected to $5$ dp and are indicated by red points in Figure 3). However as seen from Figure $3$, not all the zeros of $P_{100}(x)$ are real. This  shows that the above three conditions satisfied are not each sufficient to guarantee that all the zeros of $P_i(x)$ will be real for all $i$. 

\begin{figure}[h]\begin{center}
$ 
\begin{array}{c}
\includegraphics[height=5.5cm, width=5cm]{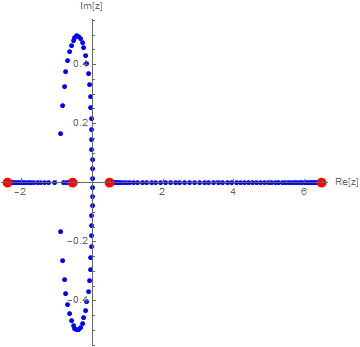}
\end{array}
$
\end{center}
\caption{The zeros of $P_{100}(x)$ for the generating function $1/(1+(-x^2+2x)t+(5x^2-1)t^2)$.}

\end{figure} 

\medskip 

\noindent
\textbf{Example 2.} Consider the sequence of polynomials $\{P_{i}(x)\}$ generated by the rational function

\begin{equation*}
 \sum_{i=0}^{\infty}P_i(x)t^i= \frac{1}{1+(2 x^2 - 8 x + 6)t+ (-5x^3+37x^2-43x-21)t^2}.
 \end{equation*}
 The corresponding $f$ is given by $$f(x)=\frac{(2 x^2 - 8 x + 6)^2}{ -5x^3+37x^2-43x-21}.$$  

\begin{figure}[h]\begin{center}
$ 
\begin{array}{c}
\includegraphics[height=5cm, width=5cm]{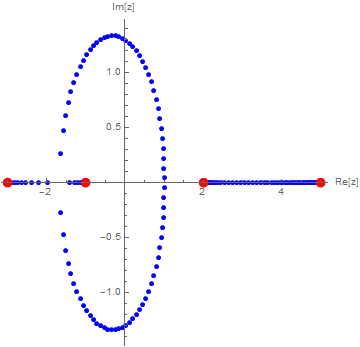}
\end{array}
$
\end{center}
\caption{The zeros of $P_{100}(x)$ for the generating function $1/(1+(2 x^2 - 8 x + 6 )t+(-5x^3+37x^2-43x-21)t^2)$.}

\end{figure} 

The zeros of $Q_1(x)$ are $1$ and $3$ which are real and simple. Also there are no ovals disjoint from $\mathbb{R}P^1$. The discriminant $D(x)= x^4-3x^3-15x^2+19x+30$ has only real zeros, i.e $x=-3, x=-1, x=2$ and $x=5.$ However $f$ has a critical value of $3.50783 \in (0,4)$ corresponding to the critical point $-1.66437$ corrected to $5$ dp, hence the condition of Theorem $5 (d)$ is violated. Consequently, some of the zeros of  $P_{100}$ are not real (see Figure $4$). The first three conditions of Theorem 5 are satisfied but not the fourth condition. Therefore having no critical value in $(0,4)$ is indeed a necessary condition for the reality of all the zeros of $P_i(x).$ \\
 
 \noindent
 \textbf{Example 3.} Consider the sequence of polynomials $\{P_{i}(x)\}$ generated by the rational function

\begin{equation*}
 \sum_{i=0}^{\infty}P_i(x)t^i= \frac{1}{1+(2 x^2 - 8 x + 6)t+ (x^4-8x^3+21x^2-14x-16)t^2}.
 \end{equation*}
 The corresponding $f$ is given by $$f(x)=\frac{(2 x^2 - 8 x + 6)^2}{ x^4-8x^3+21x^2-14x-16}.$$  

\begin{figure}[h]\begin{center}
$ 
\begin{array}{c}
\includegraphics[height=5cm, width=4.5cm]{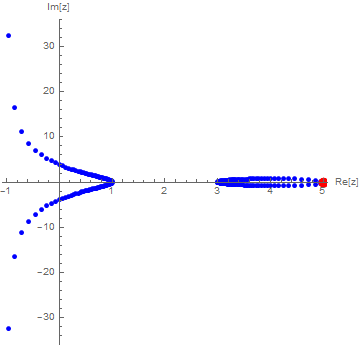}
\end{array}
$
\end{center}
\caption{The zeros of $P_{100}(x)$ for the generating function  $1/(1+(2 x^2 - 8 x + 6 )t+(x^4-8x^3+21x^2-14x-16)t^2)$.}

\end{figure} 

The zeros of $Q_1(x)$ are $1$ and $3$ which are real and simple. Also there are no ovals disjoint from $\mathbb{R}P^1$. The discriminant $D(x)= 4x^2-40x+100$ has only real zeros, i.e $x=5.$ In addition $f$ has no critical value in the real interval $(0,4)$. Thus conditions of Theorem 5 $(a)$ to $(d)$ are satisfied. Note that on the zeros of $Q_1$ we have $Q_2(1)=-16 \rlap{\kern.45em $|$} > 0$ and $Q_2(3)=-4\rlap{\kern.45em $|$} > 0$. Thus condition $(e)$ of Theorem $5$ is violated. Consequently some of the zeros of $P_{100}$ are not real as seen in Figure $5$.\\
 
 \noindent
 \textbf{Example 4.} Consider the sequence of polynomials $\{P_{i}(x)\}$ generated by the rational function

\begin{equation*}
 \sum_{i=0}^{\infty}P_i(x)t^i= \frac{1}{1+(x^2 -2x-5)t+x^2t^2}.
 \end{equation*}
The corresponding $f$ is given by $$f(x)=\frac{(x^2 -2x-5)^2}{x^2}.$$ 

\begin{figure}[h]\begin{center}
$ 
\begin{array}{c}
\includegraphics[height=5cm, width=5cm]{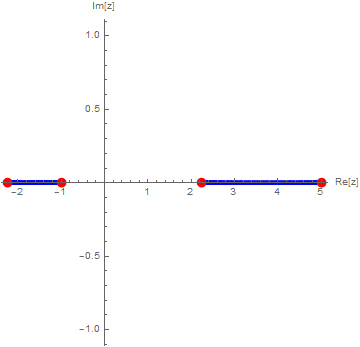}
\end{array}
$
\end{center}
\caption{The zeros of $P_{200}(x)$ when the generating function is  $1/(1+(x^2 -2x-5)t+x^2t^2)$.}

\end{figure} 

In this case, all the five conditions of Theorem $5$ have been satisfied and as seen in Figure $6$, all the zeros of $P_{200}(x)$ are real. We used $i=200,$ but any arbitrary value of $i \in \mathbb{N}^+$ works. The red points are the zeros of the discriminant and these are the endpoints of the intervals where all the zeros of $P_i(x)$ for all $i$ are located, that is, all the zeros of $P_i(x)$ for all $i$ are supported on the real axis and the support is a union of two disjoint real intervals given by $[-\sqrt{5},  -1]\cup [\sqrt{5}, 5] \subset \mathbb{R}P^1$.  The zeros of $P_i(x)$ are dense in this support as $i \to \infty.$

\section{Final Remarks}
\textbf{Problem}. Describe similar conditions guaranteeing reality of roots for all polynomials in the context of Theorem  \ref{thm:three}.

\end{document}